\documentclass[11pt]{article}
\usepackage{amssymb}
\usepackage{amsmath}
\usepackage{graphicx,color}
\usepackage{wrapfig,framed}

\newtheorem{theorem}{Theorem}
\newtheorem{guess}[theorem]{Conjecture}

\newtheorem{prop}[theorem]{Proposition}

\newtheorem{lemma}[theorem]{Lemma}

\newtheorem{remark}[theorem]{Remark}
\newtheorem{defi}[theorem]{Definition}
\newtheorem{primer}[theorem]{Example}

\newcommand{\eqa}{\begin{eqnarray}}
\newcommand{\eeqa}{\end{eqnarray}}
\newcommand{\beq}{\begin{equation}}
\newcommand{\eeq}{\end{equation}}
\newcommand{\nn}{\nonumber}

\newcommand{\F}{{\cal F}}

\newcommand{\tr}{{\rm tr}}

\newcommand{\bt}{{\bf t}}
\newcommand{\e}{\epsilon}
\newcommand{\p}{\partial}

\usepackage{array}
\newcolumntype{M}[1]{>{\centering\arraybackslash}m{#1}}
\newcolumntype{R}[1]{>{\raggedleft\arraybackslash}m{#1}}
\newcolumntype{N}{@{}m{0pt}@{}}

\def\={\,=\,}
\def\+{\,+\,}
\def\:={\,:=\,}

\begin{document}

\title{On Gromov--Witten invariants of $\mathbb{P}^1$}
\author
{
{Boris Dubrovin${}^*$, Di Yang${}^{\dag}$}\\
{\small ${}^*$ SISSA, via Bonomea 265, Trieste 34136, Italy} \\
{\small ${}^{\dag}$ Max-Planck-Institut f\"ur Mathematik, Vivatsgasse 7, Bonn 53111, Germany} \\
}
\date{}
\maketitle
\begin{abstract}
We propose a conjectural explicit formula of generating series of a new type for Gromov--Witten 
invariants of $\mathbb{P}^1$ of all degrees in \textit{full genera}.
\end{abstract}


\setcounter{equation}{0}
\setcounter{theorem}{0}
\section{Introduction}
Let $\overline{\mathcal{M}}_{g,n}(\mathbb{P}^1,\beta)$ be the moduli stack of $n$-pointed stable maps of curves of genus $g$, 
degree $\beta\in H_2(\mathbb{P}^1;\mathbb{Z})$ with target $\mathbb{P}^1$
$$
\overline{\mathcal{M}}_{g,n}(\mathbb{P}^1,\beta) \= \bigl\{\,f: \bigl(\Sigma_{g}, p_1,\dots,p_n\bigr)\rightarrow \mathbb{P}^1 \, \big| \,~ f_*\bigl([\Sigma_g]\bigr)= \beta \bigr\}/\sim.
$$ 
Here, $(\Sigma_{g},p_1,\dots,p_n)$ denotes an algebraic curve of genus $g$ with at most double-point singularities as well as with 
 the distinct marked points $p_1,\dots,p_n$, and the equivalence relation $\sim$ is defined by isomorphisms of $\Sigma_g\rightarrow \mathbb{P}^1$ 
identical on $\mathbb{P}^1$ and on the markings.  
Denote by $\mathcal{L}_i$ the $i^{th}$ tautological line bundle on $\overline{\mathcal{M}}_{g,n}(\mathbb{P}^1,\beta)$, and by 
 $\psi_i:= c_1(\mathcal{L}_i),$ $i=1,\dots,n$  the $\psi$-classes, and ${\rm ev}_i, \, i=1,\dots,n$ the evaluation maps
$$
{\rm ev}_i: \overline{\mathcal{M}}_{g,n}(\mathbb{P}^1,\beta)\rightarrow \mathbb{P}^1,  
\qquad \bigl(f: \left(\Sigma_{g}; p_1,\dots,p_n\right)\rightarrow \mathbb{P}^1\bigr)\,\mapsto f(p_i).
$$
The genus $g$, degree $\beta$ Gromov--Witten (GW) invariants
of $\mathbb{P}^1$ are integrals of the form
\beq\label{def-gw}
\int_{\left[\overline{\mathcal{M}}_{g,n}(\mathbb{P}^1,\beta)\right]^{{\rm virt}}} 
{\rm ev}_1^*(\phi_{\alpha_1})  \cdots  {\rm ev}_n^*(\phi_{\alpha_n}) \cdot \psi_1^{k_1} \cdots \psi_n^{k_n}, 
\qquad \alpha_1,\dots,\alpha_n=1,2,\, k_1,\dots,k_n\geq 0. 
\eeq
Here,  $\phi_1=1\in H^0(\mathbb{P}^1;\mathbb{C})$ is the trivial class, $\phi_2=\omega\in H^2(\mathbb{P}^1; \mathbb{C})$ is normalized by
$
\int_{\mathbb{P}^1} \omega =1
$ 
and $\left[\overline{\mathcal{M}}_{g,n}(\mathbb{P}^1,\beta)\right]^{{\rm virt}}$ denotes the virtual fundamental class \cite{KM,B,BF}.  
Clearly, the ``degree" $\beta\in H_2\bigl( \mathbb P^1; \mathbb Z\bigr)$ can be replaced by
an integer $d$ through
$d:=\int_\beta \omega$.

\noindent \textit{Notations}: ~ For any $n\geq 1$ and for a given set of integers  $k_1,\dots,k_n,\,\alpha_1,\dots,\alpha_n$, denote
\eqa
&& \hspace{-4mm} \langle \tau_{k_1}(\phi_{\alpha_1}) \dots \tau_{k_n}(\phi_{\alpha_n})  \rangle_{g,d,n} 
\:=  \int_{[\overline{\mathcal{M}}_{g,n}(\mathbb{P}^1,d)]^{\rm virt}} 
{\rm ev}_1^*(\phi_{\alpha_1})  \cdots  {\rm ev}_n^*(\phi_{\alpha_n}) \, \psi_1^{k_1} \cdots \psi_n^{k_n},  \label{corr-gd}\\
&& \hspace{-4mm} \langle \tau_{k_1}(\phi_{\alpha_1}) \dots \tau_{k_n}(\phi_{\alpha_n})  \rangle_n (\epsilon,q)
\:= \sum_{g=0}^\infty \sum_{d=0}^\infty  \epsilon^{2g-2} q^d \langle \tau_{k_1}(\phi_{\alpha_1}) \dots \tau_{k_n}(\phi_{\alpha_n})  \rangle_{g,d}. \label{corr-n}
\eeqa
Note that the sub-indices $n$ of $\langle \,,\, \rangle$ on the l.h.s. of \eqref{corr-gd} and of \eqref{corr-n} will be often omitted.

Due to the degree--dimension matching, the GW invariants \eqref{def-gw} are  \textit{zero} unless
\beq \label{dd}
2g-2+ 2d + 2 n \= \sum_{j=1}^n k_j+\sum_{j=1}^n \alpha_j.
\eeq
Hence $ \e^2\, \langle \tau_{k_1}(\phi_{\alpha_1}) \dots \tau_{k_n}(\phi_{\alpha_n})  \rangle (\epsilon,q)$ are polynomials in $\epsilon, \,q.$  More precisely, 
$$
\langle \tau_{k_1}(\phi_{\alpha_1}) \dots \tau_{k_n}(\phi_{\alpha_n})  \rangle (\epsilon,q)
\= \sum_{g+d= 1-n+\frac{k_1+\dots+k_n+\alpha_1+\dots+\alpha_n}{2}}  \epsilon^{2g-2} q^d \langle \tau_{k_1}(\phi_{\alpha_1}) \dots \tau_{k_n}(\phi_{\alpha_n})  \rangle_{g,d},
$$
where it is understood that if $k_1+\dots+k_n+\alpha_1+\dots+\alpha_n$ is an odd number then the r.h.s. vanishes.
We call $\langle \tau_{k_1}(\phi_{\alpha_1}) \dots \tau_{k_n}(\phi_{\alpha_n}) \rangle(\epsilon,q)$ the $n$-point $\mathbb{P}^1$ correlators.

In this paper, we will be particularly interested in the $\mathbb{P}^1$ correlators of the form
\beq\label{p1corrinstudy}
\left\langle  \tau_{k_1}(\omega) \dots \tau_{k_n}(\omega) \right\rangle(\e;q) \= 
\sum_{2g-2+2d=  \sum_{i=1}^n k_i}  \epsilon^{2g-2} q^d \left\langle \tau_{k_1}(\omega) \dots \tau_{k_n}(\omega)  \right\rangle_{g,d}
\eeq
(the so-called \emph{stationary sector} of the GW theory of $\mathbb P^1$ in the terminology of \cite{OP}).
These correlators vanish unless $\sum_{i=1}^n k_i $ is an even number. Due to the following quasihomogeneity 
$$
\bigl\langle  \tau_{k_1}(\omega) \dots \tau_{k_n}(\omega) \bigr\rangle\bigl(\lambda\e;\lambda^2q\bigr) \= \lambda^{\sum_{i=1}^n k_i} \left\langle  \tau_{k_1}(\omega) \dots \tau_{k_n}(\omega) \right\rangle(\e;q),\quad \forall\, \lambda\neq 0,
$$
we will often set $q=1$ in \eqref{p1corrinstudy}, and denote
for simplicity
\beq\label{simple-notation}
\bigl\langle  \tau_{k_1}(\omega) \dots \tau_{k_n}(\omega) \bigr\rangle \:= \bigl\langle  \tau_{k_1}(\omega) \dots \tau_{k_n}(\omega) \bigr\rangle(\e;1) = 
\sum_{g-1+d=\frac{k_1+\dots+k_n}{2}}  \epsilon^{2g-2} \bigl\langle \tau_{k_1}(\omega) \dots \tau_{k_n}(\omega) \bigr\rangle_{g,d}.
\eeq

\begin{defi}
Define the generating series of the $n$-point $\mathbb{P}^1$ correlators (in the stationary sector) by
\beq\label{mainsum}
C_n(\lambda_1,\dots,\lambda_n;\e) \:= \e^{n} \, \sum_{k_1,\dots,k_n\geq 0}  
\frac{(k_1+1)! \cdots (k_n+1)!}{\lambda_1^{k_1+2} \cdots \lambda_n^{k_n+2}} \bigl\langle  \tau_{k_1}(\omega) \dots \tau_{k_n}(\omega) \bigr\rangle.
\eeq
\end{defi} 

For $1$-point GW invariants (i.e. $n=1$),
the following formula has been obtained by R.\,Pandharipande \cite{P} based on the \textit{Toda conjecture}\footnote{The Toda conjecture says the partition function of GW invariants of $\mathbb{P}^1$ is a tau function of the (extended) Toda hierarchy \cite{Du1, EY, EHY, Ge, CDZ}; this conjecture has been proven in \cite{DZ-toda, OP}. }
\beq\label{pandformula}
\bigl\langle \tau_{2g-2+2d}(\omega) \bigr\rangle_{g,d} \= 
\frac{1}{d!^2} \, {\rm Coef}\Bigl( \mathcal{S}(\epsilon)^{2d-1}, \epsilon^{2g}\Bigr),\qquad \forall \, g,d\geq 0, ~g-1+d\geq 0,
\eeq
where $\mathcal{S}(\e)$ denotes the following analytic function of $\e$:
$$
\mathcal{S}(\e) \= \frac{\sinh (\e/2)}{\e/2} \= 1+ \sum_{m\geq 1} \frac{\e^{2m}}{2^{2m}  (2m+1)!} \,=:\, \sum_{m\geq 0} c_{2m} \, \e^{2m}.
$$
The formula \eqref{pandformula} was later also proved by the Gromov--Witten/Hurwitz (GW/H) correspondence in \cite{OP}.
Formula \eqref{pandformula} gives
$\bigl\langle  \tau_{2j}(\omega)\bigr\rangle
=  \sum_{g=0}^{1+j}   \frac{\epsilon^{2g-2}}{(1+j-g)!^2} \, {\rm Coef} \bigl( \mathcal{S}(\epsilon)^{2(j-g)+1}, \epsilon^{2g}\bigr)$.
So the generating series \eqref{mainsum} of $1$-point $\mathbb{P}^1$ correlators has the form
\beq C_1(\lambda;\e) 
\= \sum_{j\geq 0}  
\frac{(2j+1)!}{\lambda^{2j+2}}  \sum_{g=0}^{1+j}  \frac{\epsilon^{2g-1}}{(1+j-g)!^2} \, {\rm Coef} \Bigl(\mathcal{S}(\epsilon)^{2(j-g)+1}, \epsilon^{2g}\Bigr).
\eeq
The first several terms for $C_1(\lambda;\e)$ as given by
\beq\label{C1formula}
C_1(\lambda;\e)  \= \frac{\frac{1}{\e} -\frac{\e}{24}}{\lambda^2} \+  \frac{\frac{3}{2 \, \e} +\frac{\e}{4}+\frac{7 \, \e^3}{960}}{\lambda^4} \+ \frac{\frac{10}{3 \, \e} +\frac{15\, \e}{4} +\frac{\e^3}{16} -\frac{31 \, \e^5}{8064}}{\lambda ^6} \+ \mathcal{O}\bigl(\lambda^{-8}\bigr).
\eeq

Let us proceed to the multi-point $\mathbb{P}^1$ correlators.

\begin{guess}(Main Conjecture) \label{thm1} Define a $2\times 2$ matrix-valued series by
\beq\label{ram}
{\mathcal R}(\lambda;\e) \:= 
\begin{pmatrix} 1 & 0\\ 0 & 0\end{pmatrix}
\+ 
\begin{pmatrix} 
\alpha(\lambda;\e) & \beta(\lambda;\e)\\
\gamma(\lambda;\e) & -\alpha(\lambda;\e)
\end{pmatrix} \in {\rm Mat}\left(2, \mathbb{Q}(\e) [[\lambda^{-1}]] \right)
\eeq
where
\begin{align}
& \alpha(\lambda;\e) \= \sum_{j=0}^\infty \frac{1}{4^j \, \lambda^{2j+2}} \sum_{i=0}^j   \e^{2(j-i)}  \, \frac{1}{ i! (i+1)!} \, 
\sum_{\ell=0}^i  (-1)^\ell (2i+1-2\ell)^{2j+1} \binom{2i+1}{\ell}, \\
& \gamma(\lambda;\e) \= Q(\lambda;\e) + P(\lambda; \e),  \\
& \beta(\lambda;\e) \= Q(\lambda;\e) - P(\lambda;\e), \\
& P(\lambda;\e) \:= \sum_{j=0}^\infty \frac{1}{4^j \, \lambda^{2j+1}} \sum_{i=0}^j   \e^{2(j-i)}  \, \frac{1}{i!^2} \, 
\sum_{\ell=0}^i  (-1)^\ell (2i+1-2\ell)^{2j} \biggl[\binom{2i}{ \ell } - \binom{2i}{\ell-1}\biggr], \\
& Q(\lambda;\e) \:=  - \frac12 \sum_{j=0}^\infty \frac{1}{4^j \, \lambda^{2j+2}} \sum_{i=0}^j   \e^{2(j-i)+1}  \, \frac{2i+1}{i!^2} \, 
\sum_{\ell=0}^i  (-1)^\ell (2i+1-2\ell)^{2j} \biggl[\binom{ 2i }{ \ell } - \binom{ 2i}{ \ell-1} \biggr].
\end{align}
Then the generating series \eqref{mainsum} for the $n$-point ($n\geq 2$) GW invariants of $\mathbb{P}^1$ have the form
\eqa
&& \hspace{-4mm} C_2(\lambda_1,\lambda_2;\e)=\frac{\tr \, \left[ {\mathcal R}(\lambda_1;\e) \, {\mathcal R}(\lambda_2;\e)\right]-1}{(\lambda_1-\lambda_2)^2}, \label{twopoint} \\
&& \hspace{-4mm} C_n(\lambda_1, \dots, \lambda_n;\e)=-\frac1{n}\sum_{\sigma\in S_n} \frac{\tr \,\left[{\mathcal R}(\lambda_{\sigma_1};\e)\dots {\mathcal R}(\lambda_{\sigma_n};\e)\right]}{(\lambda_{\sigma_1}-\lambda_{\sigma_2})\dots (\lambda_{\sigma_{n-1}}-\lambda_{\sigma_n}) (\lambda_{\sigma_n}-\lambda_{\sigma_1})}, \quad n\geq 3. \label{npoint}
\eeqa 
\end{guess}

Observe that $\alpha(\lambda; \e)$ and $P(\lambda;\e)$ are even series in $\e$, while $Q(\lambda;\e)$ is an odd series in $\e$. The parity symmetry will be helpful 
for simplifying computations.
For the reader's convenience, we give the first several terms of $\mathcal{R}(\lambda;\e)$
\begin{align}
& \mathcal{R}(\lambda;\e) \= 
\begin{pmatrix} 
1 & 0 \\
 0 & 0 \\
\end{pmatrix} 
\+ \lambda^{-1} 
\begin{pmatrix}
 0 & -1 \\
 1 & 0 \\
\end{pmatrix} 
\+ \lambda^{-2} 
\begin{pmatrix} 
1 & -\frac{\e}{2} \\
 -\frac{\e}{2} & -1 \\
\end{pmatrix}
 \+ \lambda^{-3} 
\begin{pmatrix}
 0 & -\frac{\e^2}{4}-2 \\
 \frac{\e^2}{4}+2 & 0 \\
\end{pmatrix} \nn\\
&  \qquad \qquad \qquad 
\+ \lambda^{-4} 
\begin{pmatrix}
 \frac{\e^2}{4}+3 & -\frac{\e^3}{8}-3 \e \\
 -\frac{\e^3}{8}-3 \e & -\frac{\e^2}{4}-3 \\
\end{pmatrix}
\+  \mathcal{O}\bigl(\lambda^{-5}\bigr).\nn
\end{align}

\begin{remark} 
In \cite{OP} A.~Okounkov and R.~Pandharipande obtained another interesting explicit generating  
series of the multi-point relative $\mathbb{P}^1$ correlators (see the Theorem 3 of \cite{OP}).
Their generating function is labelled by a pair of partitions $\mu$, $\nu\in\mathbb Y$. 
For the particular case of $\mu=\nu=\bigl( 1^d\bigr)$ it gives the ${\mathbb P}^1$ correlators of degree $d$.
In \cite{P}, R.~Pandharipande obtained an explicit formula for GW invariants of $\mathbb{P}^1$ in degree $1$ 
based on the Toda conjecture. The Pandharipande's formula (see below in Section \ref{section3}) 
was helpful for us to verify the Main Conjecture for some particular correlators.
\end{remark}

\begin{remark}
In \cite{NS} P.~Norbury and N.~Scott considered the generating series of $\mathbb{P}^1$ 
correlators with fixed genus, and conjectured that they satisfy an explicit recursion 
of the Chekhov-Eynard--Orantin type, which was later confirmed in \cite{DOSS}.
\end{remark}

\begin{primer}
Using the Main Conjecture
we have computed some $\mathbb{P}^1$ correlators with the help of a computer program. Let us list a few of them
\begin{align}
&  \bigl\langle \tau_1(\omega)^6 \bigr\rangle \= 120 \, \e^{-2} +40 + \frac{\e^2}{2},\label{e1}\\
&  \bigl\langle \tau_2(\omega)^5 \bigr\rangle \= 36 \, \e^{-2} +\frac{2513}{24} +\frac{9745 \, \e^2}{144} 
+\frac{5435 \, \e^4}{768}+\frac{2801 \, \e^6}{82944}+ \frac{\e^8}{7962624}, \label{e2}\\
&  \bigl\langle \tau_3(\omega)^4 \bigr\rangle \= 
\frac{\e^{-2}}{2}+\frac{209}{48}+\frac{1835 \, \e^2}{192} +\frac{34807 \,  \e^4}{6912}+\frac{32053 \, \e^6}{82944}+  \frac{625 \, \e^8}{663552}, \label{e3}\\
&  \bigl\langle \tau_4(\omega)^3 \bigr\rangle  \= 
\frac{\e^{-2}}{64} +\frac{59}{384} +\frac{4217 \, \e^2}{10240} +\frac{433 \, \e^4}{1536} +\frac{443323 \, \e^6}{14745600} 
 +\frac{1261 \, \e^8}{9830400} + \frac{\e^{10}}{7077888000}, \label{e4}\\
&  \bigl\langle \tau_6(\omega)^2 \bigr\rangle \= 
\frac{\e^{-2}}{9072} +\frac{1}{648} +\frac{791 \, \e^2}{138240} +\frac{30907 \, \e^4}{5806080} 
+\frac{94537 \, \e^6}{116121600}+\frac{1781 \, \e^8}{309657600}+ \frac{\e^{10}}{104044953600}.  \label{e5}
\end{align}
The computation of these correlators based on the Main Conjecture takes less than 1 second on an ordinary computer.
It should be noted that for $n\geq 2$ the degree $d=0$ part of the $\mathbb{P}^1$ correlator of the form \eqref{p1corrinstudy} vanishes \cite{GOP}. 
So the actual highest degree in $\e$ in $\eqref{p1corrinstudy}$ is smaller than or equal to $ -2+ \sum_{i=1}^n k_i.$
More examples will be given in Section \ref{section3}. 
\end{primer}

\paragraph{Organization of the paper} 
In Section \ref{section2} we design from the Main Conjecture an algorithm suitable for computations.
In Section \ref{section3} we check the validity of the Main Conjecture in several examples, 
which also provides several new numerical values of the so-called analogues of the polygon numbers; we also  
give a few large genus asymptotics for certain GW invariants based on the Main Conjecture.
Further remarks are given in Section \ref{section4}.

\paragraph{Acknowledgements} 
One of the authors D.Y. is grateful to Youjin Zhang for his advising, and to Maxim Smirnov for helpful discussions.

\setcounter{equation}{0}
\setcounter{theorem}{0}
\section{An algorithm for computing GW invariants of $\mathbb{P}^1$}  \label{section2}
Let us design a recursive procedure for calculating GW invariants of $\mathbb{P}^1$ based on the Main Conjecture.
This algorithm was developed in \cite{DY1} for the case of GUE correlators.

\begin{defi} \label{d-DR} Fix ${\bf b}=(b_1,b_2,b_3,\dots)$ an arbitrary sequence of positive integers. 
Define recursively a family of Laurent series 
$R_{K}^{\bf b}(\lambda;\e)\in {\rm Mat} \left(2,\mathbb{Q}[\e]((\lambda^{-1}))\right)$ with $K=\{k_1,\dots,k_{m}\}$ by
\eqa
\!\!\!\!\! && R^{\bf b}_{\{\}}(\lambda;\e) \:= \mathcal{R}(\lambda;\e), \label{int-iniRR}\\
\!\!\!\!\! && R^{\bf b}_{K}(\lambda;\e) \:=  \sum_{I\sqcup J=K - \{k_1\}}  \,  
\biggl[R^{\bf b}_{I}(\lambda;\e), \Bigl(\lambda^{b_{k_1}}\,R^{\bf b}_{J}(\lambda;\e)\Bigr)_+\biggr]. \label{int-rec}
\eeqa
Here $k_1,\dots,k_{m}$ are distinct positive integers, $m=|K|$, and $\mathcal{R}(\lambda;\e)$ is defined by eq.\,\eqref{ram}.
\end{defi}
\begin{lemma} \label{l-DR}
In the particular case of $b_1=b_2=b_3=\dots=b$ we have
$$
R^{\bf b}_{K}(\lambda;\e) \= R^{\bf b}_{K'}(\lambda;\e) \,=: \, R^b_{|K|}(\lambda;\e), \qquad \mbox{as long as ~} |K|=|K'|. 
$$
Moreover, the following formulae hold true for $R^b_m(\lambda;\e),\,m\geq 1$
\beq R^b_m(\lambda;\e)  \=  \sum_{i=0}^{m-1}  
 \, \binom{m-1}{i} 
\biggl[R^b_{i}(\lambda;\e), \Bigl(\lambda^b\, R^b_{m-1-i}(\lambda;\e)\Bigr)_+\biggr].\nn
\eeq
\end{lemma}
{\it Proof.} ~
The two statements follow easily from~\eqref{int-iniRR}--\eqref{int-rec}. 
\hfill $\square$

\begin{prop}[*] \label{p-DR} Let ${\bf b}=(b_1,b_2,b_3,\dots)$ be a sequence of positive integers, and $K=\{k_1,\dots,k_m\}$ a finite set of positive integers.
The following formula holds true for GW-invariants of $\mathbb{P}^1$:
\beq
\sum_{i,j\geq 1} \bigl\langle \tau_{b_{k_1}-1}(\omega) \dots \tau_{b_{k_m}-1}(\omega) \tau_{i-1}(\omega) \tau_{j-1}(\omega) \bigr\rangle \prod_{r=1}^m b_{k_r}! \frac{i! j!}{\lambda_1^{i+1}\lambda_2^{j+1}}  \=  \sum_{I \sqcup J=K} \,
\, \frac{\tr \, R_{I}^{{\bf b}}(\lambda_1;\e)\, R_{J}^{{\bf b}}(\lambda_2;\e)}{(\lambda_1-\lambda_2)^2} - \frac{\delta_{m,0} }{(\lambda_1-\lambda_2)^2}.
\eeq
Here $m=|K|$. In the particular case that $b_1=b_2=\dots=b$ for some $b\geq 1$, we have $\forall\,m\geq 0$,
\beq\label{ij-bk}
\sum_{i,j\geq 1} \bigl\langle \tau_{b-1}(\omega)^m \, \tau_{i-1}(\omega) \tau_{j-1}(\omega) \bigr\rangle \frac{i! j! b!^m}{\lambda_1^{i+1}\lambda_2^{j+1}}  \=  \sum_{i=0}^m \,\binom{m}{i} 
\, \frac{\tr \, R^b_{i}(\lambda_1;\e)\, R^b_{m-i}(\lambda_2;\e)}{(\lambda_1-\lambda_2)^2} - \frac{\delta_{m,0} }{(\lambda_1-\lambda_2)^2}.
\eeq
\end{prop}
Here and below
a proposition marked with ``\,*\," means it is a consequence of the Main Conjecture.

\setcounter{equation}{0}
\setcounter{theorem}{0}
\section{Examples} \label{section3}
In this section, we give verifications and several applications of the conjectural formulae \eqref{twopoint}--\eqref{npoint}.
\setcounter{equation}{0}
\setcounter{theorem}{0}
\subsection{Degree $1$ GW invariants of $\mathbb{P}^1$.}
Consider the GW invariants of $\mathbb{P}^1$ of degree $d=1$ in the stationary sector:
\beq\label{degree1P1}
\langle \tau_{k_1}(\omega) \dots \tau_{k_n}(\omega) \rangle_{g,1}.
\eeq
The dimension-degree matching reads  
$2g= \sum_{i=1}^n k_i$.
It is known that \eqref{degree1P1} vanishes if any of $\{k_1,\dots,k_n\}$ is an odd number \cite{P}. 
As in the Introduction, define 
$$c_{2m}:=\frac{1}{2^{2m}  (2m+1)!}, ~ m\geq 0.$$
Pandharipande obtains \cite{P} the following interesting formula
for degree $1$ GW invariants of $\mathbb{P}^1$ 
\beq\label{pand-d1}
\langle \tau_{2j_1}(\omega) \dots \tau_{2j_n}(\omega) \rangle_{g,1}  \= \prod_{i=1}^n c_{2j_i},
\eeq
where $j_1,\dots,j_n$ are arbitrary non-negative integers satisfying
$
\sum_{i=1}^n j_i = g.
$

Note that 
$$
\langle \tau_{k_1}(\omega)\dots \tau_{k_n}(\omega) \rangle_{g,1}\quad \mbox{ with } \sum_{i=1}^n k_i= 2g
$$
is the coefficient of the $\e^{2g-2}$-term in $\langle \tau_{k_1}(\omega)\dots \tau_{k_n}(\omega) \rangle$.
Then one can easily verify the correctness of the degree $1$ invariants computed in the particular examples \eqref{e1}--\eqref{e5} by 
comparing the numbers with those computed from Pandharipande's formula.

\setcounter{equation}{0}
\setcounter{theorem}{0}
\subsection{Analogues of polygon numbers}
In this subsection we are interested in computing the GW invariants of $\mathbb{P}^1$ of the form
\beq \label{polygonnumbers}
\bigl\langle \tau_b(\omega)^n \bigr\rangle_{g,d}\,, \qquad \mbox{ with } ~ n b =2g-2+2d.
\eeq
Here, $b$ is a given non-negative integer. We call them analogues of polygon numbers \cite{DY1}.

\paragraph{The case $b=0.$}
In this case,  the numbers  $\left\langle \tau_0(\omega)^n \right\rangle_{g,d}$ are primary GW invariants of $\mathbb{P}^1$. 
We obtain from the Main Conjecture that $\forall\, n\geq 2,$
$$
\left\langle \tau_0(\omega)^n \right\rangle_{g,d} \= \left\{ \begin{array}{cc} 1 & g=0,\,d=1\\ 0 & \mbox{otherwise} \end{array} \right..
$$
This agrees with the well-known fact that higher genus primary GW invariants of $\mathbb{P}^1$ vanish.

\paragraph{The case $b=1.$}
The numbers \eqref{polygonnumbers} become quite non-trivial already for $b=1$. 
In this case, the corresponding analogues of polygon numbers coincide with the classical Hurwitz numbers. 
More precisely, let $H_{g,d}$ denote the (weighted) number of genus $g$ curves 
which are $d$-sheeted covers of $\mathbb{P}^1$ with a fixed general branch divisor. 
$H_{g,d}$ are famously known as the classical Hurwitz numbers, as they were originally 
introduced and studied by Hurwitz in the beautiful papers \cite{Hur1,Hur2}. And it was proven by R.\,Pandharipande \cite{P} 
that
$$
\Bigl\langle \tau_1(\omega)^{2g+2d-2} \Bigr\rangle_{g,d} = H_{g,d} \, .
$$
We list in Table \ref{table1} the first few classical Hurwitz numbers computed from the Main Conjecture. 
For $g\leq2$ these numbers agree with the computation in \cite{OP2}.

More examples are presented in the tables \ref{table2}--\ref{table6}.
To the best of our knowledge, most of the numbers $\langle \tau_b(\omega)^n \rangle_{g,d}$ 
presented in these tables for $g\geq 3$ and $d\geq 2$ are not available from the literature;
even for $g=0,1,2$ not many of these numbers are computed out in the literature, although 
there exist several known algorithms \cite{DZ-toda, NS}.

Looking at the numbers in these tables, 
we observe an interesting phenomenon for these rational numbers $\langle \tau_b(\omega)^n \rangle_{g,d}$: they have {\it integrality!}\,\footnote{ 
Don Zagier observed a similar phenomenon in the study of higher genus FJRW invariants \cite{BDY2, BDY3}; we are grateful to him for sharing to us his knowledge 
about integrality and his observation of integrality.}. Namely, through a direct checking, 
we observe that the denominators of these numbers always contain small prime factors only, but the numerators contain large primes; 
moreover the growth of these numbers seems to be under 
certain control. 
Deriving particularly closed formulae (it would be nice if they give rise to a polynomial time in $g,d$ algorithm) for these 
numbers for $b\geq 2$ will be extremely interesting (in the $b=1$ case a polynomial time algorithm was recently 
found in \cite{DYZ} after 125 years' discovery of these numbers by Hurwitz \cite{Hur1}).

\begin{table}[!htbp]
\begin{center}\tiny 
    \begin{tabular} {|c|M{2.0cm}|M{2.0cm}|M{2.1cm}|M{2.1cm}|M{2.0cm}|M{1.8cm}|N}
    \hline
    $n$ & $g=0$ & $g=1$ & $g=2$ &$g=3$ & $g=4$ & $g=5$ & \\[8pt]
    \hline
    $2$ & 1/2 & 0   & 0  & 0 & 0 & 0 & \\[8pt]
   \hline
    $4$ & 4 & 1/2 &  0  & 0  & 0 & 0 & \\[8pt]
    \hline
    $6$ & 120 & 40 & 1/2  & 0 & 0  &  0 & \\[8pt]
    \hline
    $8$  & 8400  & 5460 &  364 & 1/2  & 0 & 0 & \\[8pt]
    \hline
    $10$ & 1088640    & 1189440  & 206640  &  3280  & 1/2 &  0 & \\[8pt]
    \hline
    $12$ & 228191040   &  382536000  &  131670000 &  7528620  & 29524  &  1/2 & \\[8pt]
    \hline
     $14$ & 70849658880   & 171121991040 &  100557737280 &  13626893280  &  271831560  & 265720   & \\[8pt]
    \hline
       $16$ & 30641612601600   & 101797606310400 &  92919587080320 &  24109381296000  &  1379375197200   & 9793126980   & \\[8pt]
    \hline
     $18$ & 17643225600000 000   & 77793710054860 800 &  103292024327331 840 &  45097329069112 320   &  5576183206513920   & 138543794363520   & \\[16pt]
     \hline
     $20$ & 1306502906183354 8800   & 743134101959208 96000 &  136749665725094 822400 &  9213770950232808 9600    &  2084792554739198 3040   & 127011635761701 6000   & \\[16pt]
    \hline
    \end{tabular}
\end{center}
\caption{$\left\langle \tau_1(\omega)^{n} \right\rangle_{g, \,d=\frac{n}2+1-g}$.} \label{table1}
\end{table}

\begin{table}[!htbp]
\begin{center}\tiny 
    \begin{tabular} {|c|M{2.0cm}|M{2.0cm}|M{2.0cm}|M{2.0cm}|M{2.0cm}|M{2.0cm}|N}
    \hline
     $n$ & $g=0$ & $g=1$ & $g=2$ &$g=3$ & $g=4$ & $g=5$ & \\[8pt]
    \hline
    1 & $\frac14$ & $\frac1{24}$ &  $\frac{7}{5760}$   & 0 & 0    & 0 & \\[8pt]
    \hline
    2 & $\frac13$ & $\frac16$& $\frac1{576}$     & 0      & 0  & 0  & \\[8pt]
    \hline
    3 & 1  & $\frac{25}{24}$  & $\frac{19}{192}$  &$\frac1{13824}$  & 0 & 0 & \\[8pt]
    \hline
    4 & 5 & $\frac{55}6$     & $\frac{263}{96}$  & $\frac{25}{432}$ & $\frac1{331776}$ & 0 & \\[8pt]
    \hline
    5 & 36 & $\frac{2513}{24}$   & $\frac{9745}{144}$ & $\frac{5435}{768}$ & $\frac{2801}{82944}$ & $\frac1{7962624}$ & \\[8pt]
    \hline
    6 & 343 & 1474 & $\frac{328033}{192}$     & $\frac{207985}{432}$  & $\frac{ 225751 }{12288}$ & $\frac{ 817 }{41472}$ & \\[8pt]
    \hline
    7 & 4096 & $\frac{592513}{24}$  
    & $\frac{366723}8$ & $\frac{364153055 }{13824}$ & $\frac{1107239 }{324}$  & $\frac{4713415}{98304}$ & \\[8pt]
    \hline
    8 & 59049 & $\frac{1439180}3$ & $\frac{190470301}{144}$  & $\frac{2648233}{2}$ & $\frac{66481768255}{165888}$  &  $\frac{378470995 }{15552}$ & \\[8pt]
    \hline
    9 & 1000000 & $\frac{84897195}8$ & 41142049  & $\frac{74726723365}{1152}$ & $\frac{597185127}{16}$  &  $\frac{2690321702971}{442368}$ & \\[8pt]
    \hline
    \end{tabular}
\end{center}
\caption{$\langle\tau_2(\omega)^n \rangle_{g, \,d=n+1-g}$.} \label{table2}
\end{table}

\begin{table}[!htbp]
\begin{center}\tiny 
    \begin{tabular} {|c|M{2.0cm}|M{2.0cm}|M{2.0cm}|M{2.0cm}|M{2.0cm}|M{2.0cm}|N}
    \hline
     $n$ & $g=0$ & $g=1$ & $g=2$ &$g=3$ & $g=4$ & $g=5$ & \\[8pt]
    \hline
    2 & $\frac1{16}$ & $\frac1{8}$ &  $\frac{25}{1152}$   & 0 & 0    & 0 & \\[8pt]
    \hline
    4 & $\frac12$ & $\frac{209}{48}$& $\frac{1835}{192}$     & $\frac{34807}{6912}$   
    & $\frac{32053}{82944}$  & $\frac{625}{663552}$  & \\[8pt]
    \hline
    6 & $\frac{333}{16}$ & $\frac{7325}{16} $ & $ \frac{1313519}{384}$ & $\frac{46028125}{4608} $  & $\frac{1176074965}{110592}$ & $\frac{2225242915}{663552} $ & \\[8pt]
    \hline
    8 & $\frac{9065}{4}$ & $\frac{1571255}{16} $ & $\frac{320152903}{192}$ & $\frac{93077990807}{6912}$ & $\frac{215408105005}{4096}$ & $\frac{5199315506441}{55296}$ & \\[8pt]
    \hline
    10 & $ \frac{3855285 }{8}$ & $\frac{1140753285 }{32}$ & $\frac{143868323725 }{128} $ & $\frac{9601626378785 }{512}$ & $\frac{177927208378767}{1024}$ 
    & $\frac{784631685765104095}{884736}$ & \\[8pt]
    \hline
    \end{tabular}
\end{center}
\caption{$\langle\tau_3(\omega)^n \rangle_{g, \,d=\frac{3n}2+1-g} $.} \label{table3}
\end{table}

\begin{table}[!htbp]
\begin{center}\tiny 
    \begin{tabular} {|c|M{1.2cm}|M{1.4cm}|M{1.6cm}|M{1.8cm}|M{2.0cm}|M{2.0cm}|M{2.0cm}|N}
    \hline
     $n$ & $g=0$ & $g=1$ & $g=2$ &$g=3$ & $g=4$ & $g=5$ &  $g=6$ &\\[8pt]
    \hline
    1 & $\frac1{36}$ & $\frac1{32}$ &  $\frac1{1920}$   & $-\frac{31}{967680}$ & 0    & 0 & 0& \\[8pt]
    \hline
    2 & $\frac1{80}$ & $\frac{5}{96}$& $\frac{421}{11520}$     & $\frac{31}{15360}$      & $\frac{1}{3686400}$  & 0  & 0& \\[8pt]
    \hline
    3 & $\frac{1}{64}$ & $\frac{59}{384}$ & $\frac{4217 }{10240}$ & $\frac{433 }{1536}$ & $\frac{443323}{14745600}$ & $\frac{1261 }{9830400}$ & $\frac{1}{7077888000}$ &\\[8pt]
    \hline
    4 & $\frac{9}{256}$ & $ \frac{21}{32} $ & $\frac{127787}{30720} $ & $\frac{900707}{92160}$ & $\frac{14478481}{1966080} $ & 
    $\frac{311747}{245760}$ & $\frac{57610061}{2359296000}$ & \\[8pt]
    \hline
    5 &$\frac{121 }{1024} $ & $ \frac{5651}{1536} $ & $\frac{1446187 }{32768}$ & $\frac{1451959 }{6144} $ & $\frac{12797341609}{23592960}$ 
    & $\frac{5503855157}{11796480}$ & $\frac{266585680493 }{2264924160} $ & \\[8pt]
    \hline
   \end{tabular}
\end{center}
\caption{$\langle\tau_4(\omega)^n \rangle_{g, \,d=2n+1-g}$.} \label{table4}
\end{table}

\begin{table}[!htbp]
\begin{center}\tiny 
    \begin{tabular} {|c|M{1.0cm}|M{1.15cm}|M{1.4cm}|M{1.5cm}|M{2.2cm}|M{2.35cm}|M{2.55cm}|N}
    \hline
     $n$ & $g=0$ & $g=1$ & $g=2$ &$g=3$ & $g=4$ & $g=5$ & $g=6$  & \\[8pt]
    \hline
    2 & $\frac1{864}$ & $\frac1{96}$ &  $\frac{451}{23040}$   & $\frac{2597}{414720}$ &  $\frac{8281}{66355200}$  & 0 & 0& \\[8pt]
    \hline
    4 & $\frac1{1728}$ & $\frac{1039}{41472}$& $\frac{12161}{31104}$     & $\frac{8658131}{3317760}$      & $\frac{80902129}{11059200}$  & $\frac{6108849167}{796262400}$  
    & $\frac{28686913747}{11943936000}$  &  \\[8pt]
    \hline
    6 & $\frac{137 }{82944}$ & $\frac{46691 }{248832}$ & $\frac{72455425 }{7962624}$ & $\frac{3734329163 }{15925248}$ 
    & $\frac{3231504856837 }{955514880}$ & $\frac{311933225742569 }{11466178560} $ & $ \frac{108033950880129851}{917294284800}$ & \\[8pt]
    \hline
    8 & $\frac{113507}{8957952} $ & $\frac{103619845}{35831808} $ & $\frac{164491428073 }{537477120}$ & $\frac{6803735203921 }{358318080}$ 
    & $\frac{127548309823336381}{171992678400}$ & $\frac{5129142288162642911 }{275188285440}$ & $\frac{3730500946382673048971 }{12383472844800}$ & \\[8pt]
    \hline
 \end{tabular}
\end{center}
\caption{$\langle\tau_5(\omega)^n \rangle_{g, \,d=\frac{5n}2+1-g}$.} \label{table5}
\end{table}

\begin{table}[!htbp]
\begin{center}\tiny 
    \begin{tabular} {|c|M{1.1cm}|M{1.2cm}|M{1.6cm}|M{2.0cm}|M{2.0cm}|M{2.0cm}|M{2.2cm}|N}
    \hline
     $n$ & $g=0$ & $g=1$ & $g=2$ &$g=3$ & $g=4$ & $g=5$ & $g=6$  & \\[8pt]
    \hline
    1 & $\frac1{576}$ & $\frac5{864}$ &  $\frac{13}{7680}$   & $\frac{1}{322560}$ &  $\frac{127}{154828800}$  & 0 & 0& \\[8pt]
    \hline
    2 & $\frac1{9072}$ & $\frac{1}{648}$& $\frac{791}{138240}$     & $\frac{30907}{5806080}$      & $\frac{94537}{116121600}$  & $\frac{1781}{309657600}$  
    & $\frac{1}{104044953600}$  &  \\[8pt]
    \hline
    3 & $\frac1{46656}$ & $\frac{31}{41472}$     & $\frac{15431}{1658880}$  & $\frac{13082513}{278691840}$ & $\frac{55549391}{619315200}$ 
    & $\frac{114802747}{2123366400}$ & $\frac{44854036799}{6242697216000}$ & \\[8pt]
    \hline
    4 & $\frac{13}{1679616}$  & $\frac{197}{373248}$  & $\frac{1324607}{89579520}$  &$\frac{191700403}{940584960}$  
    & $\frac{62268350861}{44590694400}$ & $\frac{150956609173}{33443020800}$ & $\frac{99806823299633}{16052649984000}$ & \\[8pt]
    \hline
    5 & $\frac{1}{236196}$ & $\frac{8789}{17915904}$ & $\frac{8175239 }{322486272}$ & $\frac{19383629785 }{27088846848}$ & 
    $\frac{461054026649 }{40131624960}$ & $\frac{2400460683943939 }{23115815976960}$ & $\frac{246762110732615767 }{485432135516160}$ &  \\[8pt]
    \hline
 \end{tabular}
\end{center}
\caption{$\langle\tau_6(\omega)^n \rangle_{g, \,d=3n+1-g}$.} \label{table6}
\end{table}

\setcounter{equation}{0}
\setcounter{theorem}{0}
\subsection{Further examples. Some large genus asymptotics.}
The following proposition gives some simple consequences of the Main Conjecture.
\begin{prop}[*]
\begin{align}
& \e^2\, \sum_{k\geq 0}  \langle \tau_0(\omega) \tau_k(\omega)   \rangle   \frac{1!\, (k+1)!}{\lambda^{k+2}}  
\= \alpha(\lambda;\e) \= \frac{1}{\lambda^2}+ \frac{\frac{\e^2}{4}+3}{\lambda^4}+\frac{\frac{\e^4}{16}+\frac{15 \,\e^2}{2}+10}{\lambda^6}+\mathcal{O}(\lambda^{-8}), \label{i1simple}\\
& \e^2 \, \sum_{k\geq 0}  \langle \tau_1(\omega) \tau_k(\omega)   \rangle   \frac{2!  \, (k+1)!}{\lambda^{k+2}}  
 \= 2 \, \lambda \, \alpha(\lambda;\e) + \beta(\lambda; \e) - \gamma(\lambda;\e),\label{i2simple}\\
&
 \e^2\, \sum_{k\geq 0}  \langle \tau_2(\omega) \tau_k(\omega)   \rangle   \frac{3! \, (k+1)!}{\lambda^{k+2}}  
 \= 1 + (3 \lambda^2+2) \, \alpha(\lambda;\e) + \frac12 ( 4 \lambda - \e) \beta(\lambda;\e) - \frac12 (4 \lambda+\e) \gamma(\lambda;\e),    \label{i3simple} \\
&  \e^2\, \sum_{k\geq 0}  \langle \tau_3(\omega) \tau_k(\omega)   \rangle   \frac{4! \, (k+1)!}{\lambda^{k+2}} 
\=  2 \lambda + ( 4 \lambda^3 +4 \lambda) \, \alpha(\lambda;\e)+\frac{1}{4} \bigl(\e^2-4 \e \lambda+12 \lambda^2+8\bigr) \, \beta(\lambda;\e) \nn\\ 
&\qquad\qquad\qquad\qquad\qquad\qquad\qquad -\frac{1}{4} \bigl(\e^2+4 \e \lambda+12 \lambda^2+8\bigr) \, \gamma(\lambda;\e), \label{i4simple}\\
& \e^3\, \sum_{k\geq 0}  \langle \tau_0(\omega) \tau_1(\omega) \tau_k(\omega)   \rangle   \frac{1! \, 2! \, (k+1)!}{\lambda^{k+2}}  
\= - \frac12 (2\lambda-\e) \, \beta(\lambda;\e)- \frac12(2\lambda+\e) \, \gamma(\lambda;\e) ,   \label{i01simple} \\
&  \e^3\, \sum_{k\geq 0}  \bigl\langle \tau_1(\omega)^2  \tau_k(\omega)   \bigr\rangle   \frac{2!^2 \, (k+1)!}{\lambda^{k+2}}  
\=  \e + 2 \, \e \, \alpha(\lambda;\e) - \frac14 (2\lambda-\e)^2 \, \beta(\lambda;\e) - \frac14 (2\lambda+\e)^2 \, \gamma(\lambda;\e).   \label{i11simple}
\end{align}
\end{prop}

We want to apply these identities to study some large genus behaviour of certain $\mathbb{P}^1$ correlators.
Indeed, comparing the coefficients of both sides of \eqref{i1simple} we obtain
$$
\langle \tau_0(\omega) \tau_{2g+2d-2}(\omega) \rangle_{g, d}  =  \frac{1}{2^{2g+2d-2} \, (2g+2d-1)! (d-1)! d!} \, 
\sum_{\ell=0}^{d-1}  (-1)^\ell (2d-1-2\ell)^{2g+2d-1} \binom{ 2d-1}{  \ell  }.
$$
For example,
\begin{align}
& \langle \tau_0(\omega) \tau_{2g}(\omega) \rangle_{g, \,d=1} \= \frac{1}{4^g \cdot (2g+1)! \cdot 0! \cdot 1!}, \nn\\
& \langle \tau_0(\omega) \tau_{2g+2}(\omega) \rangle_{g, \,d=2} \= \frac{3^{2g+3}-3}{4^{g+1} \cdot (2g+3)! \cdot 1!\cdot 2!}, \nn\\
& \langle \tau_0(\omega) \tau_{2g+4}(\omega) \rangle_{g, \,d=3} \= \frac{5^{2g+5} -3^{2g+5} \cdot 5 + 10}{4^{g+2} \cdot (2g+5)! \cdot 2! \cdot 3!}.\nn
\end{align}

\begin{prop}[*] For fixed $d\geq 1$, the following asymptotic holds true 
\beq
(2g+2d-1)! \, \langle \tau_0(\omega) \tau_{2g+2d-2}(\omega) \rangle_{g, \,d} \, \sim \, 2 \, \frac{(d-1/2)^{2d-1} }{d! \cdot (d-1)!} (d-1/2)^{2g},\quad g\rightarrow \infty.
\eeq
\end{prop}

Comparing the coefficients of both sides of \eqref{i2simple} we obtain
for $0\leq g\leq j$ that 
\eqa
&& \langle \tau_1(\omega) \tau_{2j-1}(\omega) \rangle_{g, \,d=1+j-g} \nn\\
&=& \frac{1}{4^j \, (2j)!}  \biggl[ \frac1{ (j-g)! (j-g+1)!} \, 
\sum_{\ell=0}^{j-g}  (-1)^\ell (2(j-g)+1-2\ell)^{2j+1}   \binom{ 2(j-g) +1}{\ell} \nn\\
&& \qquad  -\frac1{(j-g)!^2} \sum_{\ell=0}^{j-g}  (-1)^\ell (2(j-g)+1-2\ell)^{2j}  \biggl( \binom{ 2(j-g) }{\ell} - \binom{2(j-g)}{\ell-1}\biggr) \biggr].\nn
\eeqa
In other words, for $g\geq0, \, d\geq 1$, we have 
\eqa
&&\langle \tau_1(\omega) \tau_{2g+2d-3}(\omega) \rangle_{g, \,d} \nn\\
&=& \frac{1}{4^{g+d-1} \, (2g+2d-2)!}  \biggl[ \frac1{ (d-1)! \, d!} \, 
\sum_{\ell=0}^{d-1}  (-1)^\ell (2d-1-2\ell)^{2g+2d-1}   \binom{ 2d-1}{\ell} \nn\\
&& \qquad  -\frac1{(d-1)!^2} \sum_{\ell=0}^{d-1}  (-1)^\ell (2d-1-2\ell)^{2g+2d-2}  \biggl( \binom{ 2d-2 }{\ell} - \binom{2d-2}{\ell-1}\biggr) \biggr].\nn
\eeqa
Clearly, if $d=1$ then the above formula gives  $\langle \tau_1(\omega) \tau_{2g-1}(\omega) \rangle_{g, \,1}\equiv0$. We arrive at

\begin{prop}[*] For fixed $d\geq 2$, the following asymptotic holds true 
\beq
(2g+2d-2)! \, \langle \tau_1(\omega) \tau_{2g+2d-3}(\omega) \rangle_{g, \,d} \;\sim\; (d-1)\frac{(d-1/2)^{2d-2}}{d! \, (d-1)!} (d-1/2)^{2g},\quad g\rightarrow \infty.
\eeq
\end{prop}

Similarly from \eqref{i3simple} we obtain
\begin{prop}[*] For fixed $d\geq 1$, the following asymptotic holds true 
\beq
(2g+2d-3)! \, \langle \tau_2(\omega) \tau_{2g+2d-4}(\omega) \rangle_{g, \,d} \;\sim\; \frac{4d^2-6d+3}{12} \frac{(d-1/2)^{2d-3}}{d! \, (d-1)!} (d-1/2)^{2g},\quad g\rightarrow \infty.
\eeq
\end{prop}

More generally, we have
\begin{prop}[*] For fixed $ k\geq0$ and fixed $d\geq 1$,  the following asymptotic holds true 
\eqa
&& (2g+2d-k-1)! \, \langle \tau_k(\omega) \tau_{2g+2d-k-2}(\omega) \rangle_{g, \,d} \nn\\
&& \;\sim\;  2 \frac{(d-1/2)^{2d}}{(k+1)!\, d!^2} \left(1+ \frac{(-1)^k}{2^{k+1} (d-1/2)^{k+1}} \right) (d-1/2)^{2g},\quad g\rightarrow \infty.
\eeqa
\end{prop}


\setcounter{equation}{0}
\setcounter{theorem}{0}
\section{Concluding remarks} \label{section4}
The first remark is on the motivation of our Main Conjecture and on an idea of a possible proof. 

\paragraph{Toda Conjecture.}{\bf A. Weak version}: \textit{Let $\F$ be the following generating series of 
GW-invariants of $\mathbb{P}^1$  (often called the free energy)
\beq
\F=\F(x,\bt; \e):=\sum_{g,d\geq 0} \e^{2g-2} \sum_{m,n_0,n_1,n_2,\dots\geq 0}  \frac{x^m}{m!} \prod_{j=0}^\infty \frac{t_j}{n_j!} \,
\Bigl\langle \tau_0(1)^m \prod_{j=0}^\infty \tau_j(\omega)^{n_j} \Bigr\rangle_{g,d} \,.
\eeq
Here $\bt=(t_0,t_1,t_2,\dots).$ Let $Z:=e^{\F}$. Define $u,v$ by
\eqa
v =v(x,\bt;\e) & := & \e \frac{\p }{\p t_0 } \log \frac{Z(x+\e, \bt; \e)}{Z(x,\bt;\e)}, \nn\\
u = u(x,\bt;\e) & := & \log \frac{Z(x+\e, \bt; \e) \, Z (x-\e, \bt;\e)}{Z^2(x,\bt;\e)}. \nn
\eeqa
Then $u,v$ satisfy the Toda hierarchy with the first equation being
\eqa
\frac{\p v(x,\bt;\e)}{\p t_0} &=& \frac1{\e} \, \Bigl(e^{u(x+\e,\bt;\e)} - e^{u(x,\bt;\e)}\Bigr),\nn\\
\frac{\p u(x,\bt;\e)}{\p t_0} &=& \frac1{\e} \, \bigl(v(x,\bt;\e)-v(x-\e,\bt;\e)\bigr).\nn 
\eeqa
{\bf B. Strong version}: $Z$ is a tau function (in the sense of \cite{DZ-toda,DY1}) of the Toda hierarchy.}
Note that the strong version of the Toda conjecture along with the following celebrated string equation 
\beq\label{stringeq}
\sum_{i=1}^\infty t_i \frac{\p Z}{\p t_{i-1}}  \+ \frac{x \, t_0}{\e^2} \= \frac{\p Z}{\p x}
\eeq
uniquely determines $Z$ up to only a constant factor (independent of $\e$\,!). 
Validity of the Toda conjecture was confirmed in \cite{OP, OPeq} (see also \cite{DZ-toda}).
We would like to mention 
that an extended version of the Toda conjecture, which contains the full information of GW invariants 
of $\mathbb{P}^1$, was obtained in \cite{CDZ}. At present, we do not know how to generalize our Main Conjecture
to the extended Toda hierarchy of \cite{CDZ}. Let us continue to explain our motivation. 
In \cite{DY1} we derived explicit generating series in terms of matrix resolvents for 
logarithmic derivatives of tau-function of \textit{an arbitrary solution} to 
the Toda hierarchy. Our simple but main observation that motivates the Main Conjecture of the present paper is that 
the particular solution corresponding to GW invariants of $\mathbb{P}^1$ (the Toda conjecture) is characterized by the following initial data
\eqa
u(x,\bt= {\bf 0}; \e) &=& 0\,, \label{ini1}\\
v(x,\bt= {\bf 0}; \e) &=& x+ \frac \e 2\,. \label{ini2}
\eeqa
This can be deduced from the string equation combined with the divisor equation for GW invariants. 
Hence one can expect to prove the Main Conjecture  
by using the results of \cite{DY1} about matrix resolvents of difference operators. Indeed, using \eqref{ini1}--\eqref{ini2} 
and the results of \cite{DY1}, it is straightforward to reduce the computation of the GW invariants of $\mathbb{P}^1$ 
in the stationary sector to the following
problem of finding the unique matrix-valued formal series $R_n(\lambda)$:
\begin{align}
& R_{n+1}(\lambda) \, U_n(\lambda) \,-\, U_n(\lambda) R_n(\lambda) \= 0\,, \label{r1} \\
& {\rm tr} \, R_n(\lambda) \= 1,\qquad {\rm \det} \, R_n(\lambda) \= 0\,,  \label{r2} \\
& R_n(\lambda) \= \begin{pmatrix} 1 & 0  \\  0 & 0 \\ \end{pmatrix}
\+ \mathcal{O}\bigl(\lambda^{-1}\bigr) \, \in \, {\rm Mat} \bigl(2, \mathbb{Z}[n][[\lambda^{-1}]]\bigr)  \label{r3}
\end{align}
where $U_n(\lambda) := \begin{pmatrix} n\e +\frac\e2-\lambda & 1\\ -1 & 0\end{pmatrix}$.
Recently, in \cite{DYZ2} we achieve a proof of the Main Conjecture by using this idea.

In \cite{Du1,DZ,Du2} the first-named author of the present paper and Y.~Zhang developed an approach of 
computing GW invariants of any smooth projective variety with semisimple quantum cohomology. Recall that 
the quantum cohomology of $\mathbb{P}^1$ is semisimple, hence one can apply the Dubrovin--Zhang approach 
to compute in principle all the GW invariants \eqref{def-gw} of $\mathbb{P}^1$; see in \cite{DZ-toda} for the details. 
Our examples (see Section \ref{section3}) for $g\leq2$ based on the Main Conjecture agree with \cite{DZ-toda}.  
In particular, in \cite{DYZ} D.~Zagier and the authors of the present paper have designed several (new) 
algorithms of computing Hurwitz numbers $H_{g,d}$ (one of which is based on
the Dubrovin--Zhang approach);  one can easily verify that Table \ref{table1} agrees with the computation in \cite{DYZ}.

We believe that the method developed in \cite{BDY1,BDY2,BDY3,DY1} can be applied to and it would be very useful for the computation of
GW invariants of $\mathbb{P}^1$-orbifolds \cite{R,C}, for which we plan to do in a subsequent publication.

~~

\noindent{{\it Note added}}: Recently, O.~Marchal \cite{Mar} gave a proof of the Main Conjecture of the present paper
 by using the Chekhov--Eynard--Orantin topological recursion \cite{NS,DOSS,DMNPS}.

\end{document}